\input amstex 
\documentstyle{amsppt}
\magnification=1200

\topmatter  

\title Periodic cyclic homology as sheaf cohomology
\endtitle
\author by Guillermo Corti\~ nas*\endauthor
\affil Departamento de Matem\'atica\\
Facultad de Cs. Exactas y Naturales\\
    Universidad de Buenos Aires\\
       Argentina\endaffil
\thanks (*) Partially supported by grant BID802/OC-AR-PICT 2260 of
the Agencia Nacional de Promoci\'on Cient\'\i fica de La Argentina
and by grants UBACyT EX014 and TW79 of the Universidad de Buenos Aires.
\endthanks
\address Departamento de
Matem\' atica, 
Ciudad Universitaria Pabell\'on 1,
(1428) Buenos Aires, Argentina.
\endaddress
\email gcorti\@dm.uba.ar\endemail
\leftheadtext{Guillermo Corti\~ nas}
\rightheadtext{Periodic cyclic homology as sheaf cohomology}

\endtopmatter  
\NoBlackBoxes
\document
\redefine\lim{\operatornamewithlimits{lim}}
\define\limo{\underset{m}\to{\lim}}
\define\limu{\underset{n}\to{\lim}}

\define\fib{\twoheadrightarrow}

\define\iso{\tilde\rightarrow}

\define\holi{\operatornamewithlimits{holim}}

\define\rat{\Bbb Q}
\define\zal{\Bbb Z}
\define\faz{\partial}
\define\Alg{\operatorname{Alg}}
\define\Comm{\operatorname{Comm}}
\define\Mod{\operatorname{Mod}}
\define\Hy{\Bbb H}
\define\bgp{\zal_\infty\Cal NGl}
\document
\subhead{0. Introduction}\endsubhead This paper continues the study
of the noncommutative infinitesimal cohomology we introduced in [3]. 
This is the cohomology of sheaves on a noncommutative version of 
the commutative infinitesimal site of Grothendieck ([8]). 
Grothendieck showed that, for a smooth scheme $X$ of characteristic zero,
the cohomology of the structure sheaf on the infinitesimal site
gives de Rham cohomology:
$$
H_{inf}^*(X,O)\cong H_{dR}^*(X)\tag{1}
$$
Here we prove that, for any associative, not necessarily unital algebra 
$A$ over $\rat$, the
noncommutative infinitesimal cohomology of the structure sheaf 
modulo commutators gives periodic cyclic homology (Theorem 3.0):
$$
H_{inf}^*(A,O/[O,O])\cong HP_*(A)\tag{2}
$$
We view \thetag{2} as a noncommutative affine version of \thetag{1}. 
Indeed, the noncommutative analogue of smoothness is quasi-freedom,
and for quasi-free algebras $HP_*$ agrees with (Karoubi's definition of)
noncommutative de Rham cohomology, i.e. we have:
$$
H_{inf}^*(R,O/[O,O])\cong H_{dR}^*(R):=
H^*(\Omega(R)/[\Omega(R),\Omega(R)],d)\tag{3}
$$
for $R$ quasi-free. Here $\Omega(R)$ is the DGA of noncommutative forms.
Grothendieck's theorem includes also a description of the cohomology of the 
de Rham complex $(E\otimes_{O_X}\Omega_X,\nabla)$
associated with a bundle $E$ with a flat connection as the cohomology of
a certain module on the infinitesimal site; in fact \thetag{1} is
a obtained from this by setting $E=O_X$, $\nabla=d$. A noncommutative 
version of Grothendieck's theorem for the de Rham cohomology of flat bundles
is proved in Theorem 2.3. It expresses the cohomology of the complex
$H^*(E\otimes_{\tilde{R}}\Omega(R),\nabla)$ associated with a right 
module over the augmented algebra 
$\tilde{R}$ with a flat connection $\nabla$ as sheaf cohomology; for
$(E,\nabla)=(R,d)$ it gives:
$$
H_{inf}^*(R,O)\cong H^*(\Omega(R),d)=0
$$
We also compute infinitesimal
(hyper)cohomology with coefficients in algebraic $K$-theory. We show
that, for any $\rat$-algebra $A$, we have a natural isomorphism:
$$
H^n_{inf}(A,K_1)\cong HP_n(A)\qquad (n\ge 2) \tag{4}
$$
In section 4, we use the case $n=2$ of \thetag{4} to give a natural,
sheaf theoretic 
construction for the Chern-Connes-Karoubi character:
$$
ch_0:K_0(A)@>>>HP_0(A)\cong H^2(A,K_1)\tag{5}
$$
The cohomology groups of $K_1$ not listed in \thetag{4} appear
in an exact sequence:
$$
K_2(A)@>c_2>>HN_2(A)@>>>H_{inf}^0(A,K_1)@>>>K_1(A)@>c_1>>HN_1(A)@>>>H^1(A,K_1)
@>>>0\tag{6}
$$
The map $c$ above is the Jones-Goodwillie character to negative cyclic
homology.
The sequence \thetag{6} is obtained as the lower terms of the long
exact sequence of homotopy groups of a homotopy fibration sequence (Th.
5.0):
$$
\Bbb H_{inf}(A,\Cal K)@>>>\Cal K(A)@>c>>SCN_{\ge 1}(A)\tag{7}
$$
Here $\Cal K$ is a (functorial, fibrant) simplicial version of the 
$K$-theory space $BGl^+$ and $SCN_{\ge 1}$ is the simplicial set
the Dold-Kan correspondence associates with the negative cyclic complex
truncated below dimension $1$. 
\bigskip
As stated above, this paper is a continuation of [3]. In op.cit. the case
$*=0$ of \thetag{2} was proved, and the case $*=1$ was conjectured.
It was further conjectured that the higher infinitesimal
cohomology groups with coefficients in $O/[O,O]$ and $K_1$ 
vanished; \thetag{3} and \thetag{4} disprove both conjectures.
In [3] the sequence \thetag{6} was shown to be rationally exact
everywhere except at $H^1_{inf}(A,K_1)$. The fibration sequence
\thetag{7} was not proved in [3] --not even rationally-- and is
therefore new.
\bigskip
The rest of this paper is organized as follows. In section 1 we 
give the basic definitions and notations. Section 2 concerns
flat connections; no assumption on the characteristic of the ground
field $k$ is made in this section. We
show that for a quasi-free algebra $R$,
the category of pairs 
$(E,\nabla)$ of a right $\tilde{R}$-module and a flat 
connection is equivalent to a subcategory of
the category of modules on the infinitesimal site (Proposition 2.1)
and that the sheaf cohomology of the infinitesimal module corresponding to
$(E,\nabla)$ is the cohomology of the complex 
$(E\otimes_{\tilde{R}}\Omega(R),\nabla)$ (Theorem 2.4). The isomorphism
\thetag{2} is proved in section 3. Section 4 is devoted to the
sheaf theoretic interpretation of the Chern character \thetag{5}.
The isomorphism \thetag{4} is proved in Proposition 4.0.
In 4.1 we construct a character $ch^{inf}:K_0(A)@>>>H^2_{inf}(A,K_1)$
for algebras $A$ over an arbitrary field $k$. In theorem 4.4. we prove
that in characteristic zero, the character of 4.1 agrees with that
of Connes-Karoubi up to a normalization factor of $2$.
The fibration sequence \thetag{7} and the exact sequence \thetag{6}
are the subject of section 5 (Theorem 5.0 and Corollary 5.1.)
\bigskip
\bigskip
\subhead{1. Notations and definitions}\endsubhead
\smallskip
\subsubhead{Basic notations and conventions 1.0}\endsubsubhead
We fix a field $k$ which we call the ground field. We make no
assumptions on the characteristic of $k$, except as stated. All
tensor products are over $k$ except as indicated. Algebras
are nonunital associative $k$-algebras. We use the letters $A,B\dots$ 
for general algebras and the letters $R,S\dots$ for quasi-free algebras.
We write $\tilde{A}$ for the augmented unital algebra $A\oplus k$.
Ideals are always two-sided. Modules over unital algebras are right
modules and are always assumed unital (i.e. $1\cdot m=m$ for $m\in M$).
Curly brackets $\{\}$ are used to denote pro-objects. Thus $\{A_n\}$
is a pro-algebra having $A_n$ at level $n\in\Bbb N$ where the 
structure maps $A_n@>>>A_{n-1}$ are understood from the context.
\bigskip
\subsubhead{Cosimplicial cylinder 1.1}\endsubsubhead
We write $\Delta$ for the category of the finite ordered sets 
$[n]=\{0<\dots<n\}$ and nondecreasing monotone
maps. Recall a cosimplicial object in a category $C$ is a functor 
$\Delta\mapsto C$. Any object $A$ in a category $C$ with coproducts 
gives rise to a cosimplicial object $\Delta@>>>C$, by 
$[n]\mapsto \coprod_{i=0}^nA$; we call
this the {\it coproduct cosimplicial object} associated to $A$.
We write $*$ for the coproduct in the category of algebras and 
$Q^*A:[n]\mapsto Q^n(A):=A*\dots *A$ ($n+1$ factors) for the 
coproduct cosimplicial algebra. We write $\partial^i$ and $\mu^i$
for the cofaces and codegeneracies
of this and of most cosimplicial objects appearing in this paper.
The maps $\mu^0\mu^1\dots\mu^n:Q^n(A)@>>>A$ ($n\ge 0$) 
assemble to give a cosimplicial
map $\mu$ to the constant cosimplicial algebra. We write $q^*(A)$ 
for its kernel and set:
$$
Cyl^*(A):=Q^*(A)/q^*(A)^\infty
$$
Here the exponent $^\infty$ has the same meaning as in [7]; thus 
$Cyl^*(A)$ is a cosimplicial pro-algebra or more precisely an inverse 
system of cosimplicial algebras having 
$Cyl^*(A)_n:=Q^*A/q^*(A)^{n}$ at level $n$.
We call the pro-algebra $Cyl^1(A)$ the {\it cylinder} of $A$ and  
$Cyl^*(A)$ the {\it cosimplicial cylinder}. The name comes from the role
these pro-algebras play in the notion of nil-homotopy developed in [2]
and [3]. Roughly one thinks of $Cyl^1(A)$ as two copies of $A$ with
the `contractible space' $qA/qA^\infty$ stuck in between. 
Next we shall show that for a quasi-free algebra $R$, $Cyl^*(R)$ can
be described in terms of --noncommutative-- differential forms.
By the tubular neighborhood theorem ([4, Th. 2]; see also [4, 8.7]
and [2, Theorem 2.1]), we have an isomorphism of inverse systems
$\alpha:Cyl^1(R)_{n}=Q^1(R)/q^1(R)^{n}\iso\Omega(R)/\Omega^+(R)^{n}
\cong\oplus_{i=0}^{n-1}\Omega^i(R)$, where $\Omega^+(R)$ is the ideal 
generated by $\Omega^1(R)$.
The map $\alpha$ is natural on pairs $(R,\phi)$ of a quasi-free 
algebra and a $1$-cochain $\phi:R@>>>\Omega^2(R)$ satisfying 
$\delta(\phi)=-d\cup d$, and is such that $\alpha\faz^1(a)=a$ and 
$\alpha\faz^0(a)=a+da+\phi(a)+\dots$
($a\in R$). This isomorphism generalizes to higher codimension as 
follows. Write $W$ for the coproduct cosimplicial vectorspace on $k$;
thus $W^m=\oplus_{i=0}^{m}ke_i$, and cofaces and codegeneracies are as
follows: $\faz^i(e_j)=e_j$ if $j<i$ 
and $\faz^i(e_j)=e_{j+1}$ if $j\ge i$, $\mu^i(e_j)=e_j$ if $j\le i$ and
$\mu^i(e_j)=e_{j-1}$ if $j>i$. We write $V$ for the kernel
of the canonical map to the constant cosimplicial vectorspace $k$. Thus
$V^m=\oplus_{i=1}^mkv_i$ where $v_i=e_i-e_{i-1}$ ($i=1,\dots,m$).
We have an isomorphism:
$$
\alpha^m:Cyl^m(R)_n\iso\oplus_{i=0}^{n-1}\Omega^i(R)\otimes 
T^i(V^m)\tag{8}
$$ 
The right hand side of \thetag{8} has the algebra structure 
of the quotient of the unaugmented tensor algebra 
$\Omega_m(R):=\ker(T_{\tilde{R}}(\oplus_{i=1}^n\Omega^1(R))@>>>k)$ by the
$n$-th power of the ideal generated by 
$\oplus_{i=1}^m\Omega^1(R)=$ $\Omega^1(R)\otimes V^m$. The isomorphism
$\alpha^0$ is the identity map, and $\alpha^1$
is the map $\alpha$ defined above composed with the canonical
identification induced by $\Omega^i(R)\otimes T^i(V^1)=\Omega^i(R)$.
To explain the isomorphism $\alpha^m$ for higher $m$, we need some more
notation. Let us use the provisional notation ${Cyl'}^m(R)_n$ for the right
hand side of \thetag{8}; later on we shall identify both sides and
eliminate the 's from our notation. For $m\ge 1$ and $1\le i\le m+1$,
let $\faz'_i:{Cyl'}^m(R)@>>>{Cyl'}^{m+1}(R)$, $\faz'_i(\omega\otimes x)=
\omega\otimes{\faz^i}^{\otimes r}(x)$, 
($\omega\otimes x\in\Omega^r(R)\otimes T^r(V^m)$). Also let 
$\faz'_i:R@>>>{Cyl'}^1(R)$ $\faz'_i=\alpha^1\faz^i$, $i=0,1$. 
For $m\ge 2$, $\alpha^m$ is the map 
determined by $\alpha^m\faz^m\dots\faz^{i+1}\faz^{i-1}\dots\faz^0=
\faz'_m\dots\faz'_{i+1}\faz'_{i-1}\dots\faz'_0$ ($i=0,\dots m$). 
The same argument as in the proof of [2, Theorem 2.1] shows that $\alpha^m$
is an isomorphism for all $m\ge 0$. Thus we may --and do-- identify 
$Cyl^m(R)$ with ${Cyl'}^m(R)$ and the coface and codegeneracy maps with
their images through $\alpha^*$. For $m\ge 1$ and $1\le i\le m+1$ and
for $m=0$ and $i=0,1$, the $i$-th coface map is identified with the
map $\faz'_i$ defined above. 
For $m\ge 1$ the 
coface map $\faz_0:Cyl^m(R)@>>>Cyl^{m+1}(R)$ 
can be written as follows:
$$
\faz^0=\sum_{n=0}^\infty P_n:=
\{\sum_{n=0}^rP_n:Cyl^m(R)_{r}@>>>Cyl^{m+1}(R)_{r}\}
$$
Here $P_n(\Omega^r(R)\otimes T^r(V^*))\subset\Omega^{r+n}(R)\otimes 
T^{r+n}(V^{*+1})$.
The $0$-th component is simply the map  
$P_0(\omega\otimes x)=\omega\otimes\faz^0(x)$; for $\omega=a_0da_1\dots 
da_p\in\Omega^p(R)$ and $x_1\dots x_p\in T^p(V^m)$, we have:
$$
\multline
P_1(\omega\otimes x)=d\omega\otimes v_1\faz^1(x)
+\sum_{i=1}^p a_0da_1\dots da_{i-1}\phi(a_i) da_{i+1}\dots 
da_p\otimes\\
\faz^0(x_1)\dots\faz^0(x_{i-1})
(\faz^1(x_i)v_1+v_1\faz^1(x_i)-2v_1^2)\faz^0(x_{i+1})\dots\faz^0(x_p)
\endmultline\tag{9}
$$
The formulae above are derived as follows. By definition of the canonical
isomorphism $\alpha^*$ we have 
$a+da\otimes v_i+\phi(a)\otimes v_i^2+\dots=\faz^n\dots\faz^{i+1}\faz^{i-1}\dots\faz^0(a)$
($i\ge 1$);
putting this together with the formula for $\alpha\faz_0$ and the
cosimplicial identity 
$\faz^0\faz^n\dots\faz^{i+1}\faz^{i-1}\dots\faz^0=\faz^{n+1}\dots\faz^{i+2}\faz^{i}\dots\faz^0$
we obtain:
$$
\split
\faz^0(a+da\otimes v_i+\phi(a)\otimes v_i^2+\dots)=&
a+da\otimes v_{i+1}+\phi(a)\otimes v_{i+1}^2\dots\quad(i\ge 1)\\
\faz_0(a)=& a+da\otimes v_1+\phi(a)\otimes v_1^2+\dots
\endsplit
\tag{10}
$$
The formula for $P_0$ is immediate from these identities; the formula for 
$P_1$ is immediate from that for $P_0$ and \thetag{10}. It is possible 
to iterate this process
to obtain a general formula for $P_m$ in terms of the homogeneous parts of
degrees $\le m$ of $\faz_0:R@>>>Cyl^1(R)$, but we shall have no ocassion
for this. The description of $Cyl^*(R)$ is completed as 
one checks that the codegeneracy
maps are induced by those of the cosimplicial vectorspace $V$; thus
$\mu^i(\omega\otimes x)=\omega\otimes{\mu^i}^{\otimes r}(x)$ 
$(\omega\otimes x\in\Omega^r(R)\otimes T^r(V^*))$.
\bigskip
\subsubhead{Adic filtration 1.2}\endsubsubhead
Let $R$ be a quasi-free algebra, $I\subset R$ an ideal,
$A=R/I$, and $\pi:R\fib A$ the projection map. We also write $\pi$ for
the map $Q^*R@>>>A$ to the constant cosimplicial algebra defined as
either of the equal composites $Q^*R@>\mu>>R@>\pi>> A$,
$Q^*R@>Q^*(\pi)>>Q^*A@>\mu>>A$. Put $q^*(R,I):=\ker\pi$, 
$Cyl^*(R,I):=Q^*(R)/q^*(R,I)^\infty$ and $Cyl^*(R,I)_n:=Q^*(R)/q^*(R,I)^{n}$.
We call this the {\it relative cosimplicial cylinder}. It is easy to see 
that the isomorphism $\alpha$ defined in 1.1 induces
$Cyl^m(R,I)_n\cong\Omega_m(R)/<I+\Omega_m^+(R)>^{n}$. For $m=1$
we get $Cyl^1(R,I)_n=\Omega(R)/\Cal F^n$ where
$$
\Cal F^n=\Cal F_I^n:=<I+\Omega^+(R)>^n=
\oplus_{i=0}^\infty F^{n-i}\Omega^i(R)
$$
Here $F^p\Omega^q(R)=F_I^p\Omega^q(R):=\Omega^q(R)$ if $p\le 0$ and is 
$\sum I^{i_0}dA I^{i_1}\dots I^{i_{q-1}}dAI^{i_q}$ for $p\ge 1$,
with the sum taken over all multi-indices $i=(i_0,\dots ,i_q)$ such that
$|i|:=\sum i_j\ge p$.
For all $m\ge 0$, $n\ge 1$, we have:
$$
Cyl^m(R,I)_n=\oplus_{i=0}^{n-1}\frac{\Omega^i(R)}{F^{n-i}\Omega^i(R)}
\otimes T^i(V^m)\tag{11}
$$
For future use we also record here the behavior of the filtration $\Cal F$
with respect to the Karoubi operator and the Connes, de Rham and Hochschild
boundary maps:
$$
\kappa\Cal F^n\subset \Cal F^n \qquad B\Cal F^n\subset\Cal F^n\qquad
d\Cal F^n\subset\Cal F^{n}\qquad b\Cal F^n\subset\Cal F^{n-1}
\tag{12}
$$
\bigskip
\subsubhead{Noncommutative infinitesimal cohomology 1.3}\endsubsubhead
Let $\Cal G$ be a covariant functor going from the category of $k$-algebras
to the category of abelian groups, $R$ a quasi-free algebra, $I\subset R$ an ideal,
and $A=R/I$. Consider the cosimplicial abelian group 
$\Cal G Cyl^*(R,I)${\^\ \ }$\negthickspace:[n]\mapsto\limo\Cal GCyl^n(R,I)_m$.
We write $H_{inf}^*(A,\Cal G)$
for the cohomology of $\Cal G Cyl^*(R,I)$\^\ \ ; we call this the 
{\it infinitesimal cohomology}
of $A$ with values in $\Cal G$. 
Thus $H_{inf}^*(A,\Cal G)$
can be computed as the cohomology of the cochain complex:
$$
C(R,I,\Cal G):\lim_{n}\Cal G(R/I^n)@>\faz>>\lim_n\Cal G(Cyl^1(R,I)_n)@>\faz>>
\lim_n\Cal G(Cyl^2(R,I)_n)@>\faz>>\dots\tag{13}
$$
where $\faz=\sum_{i=0}^{m+1}(-1)^i\faz^i:C^m=C^m(R,I,\Cal G)@>>>C^{m+1}$
or equivalently --by the Dold-Kan correspondence-- as the cohomlogy of the 
normalized complex $\bar{C}^m=C^m/\sum_{i=1}^m\faz^i(C^{m-1})$
with as coboundary map the map induced --and still denoted by-- $\faz^0$.
It is implicit in our notation that the cohomology of the complex \thetag{13}
depends only on $A$ rather than on its presentation as a quotient of a 
quasi-free algebra. In fact the groups $H^*(A,\Cal G)$ have the following 
interpretation
in terms of sheaf cohomology. Consider the category $inf(\Alg/A)$ having as objects the
surjective algebra maps $p:B\fib A$ with nilpotent kernel and as morphisms
$(p:B\fib A)@>>>(p':B'\fib A)$ the homomorphisms $B@>>>B'$ making the obvious
diagram commute. Regard $inf(\Alg/A)^{op}$ as a site with the indiscrete 
Grothendieck topology; i.e. the coverings are the families $\{ B\cong B'\}$ consisting
of a single isomorphism; this site is the {\it noncommutative infinitesimal
site} over $A$. Any covariant functor
$inf(\Alg/A)@>>>Ab$ to abelian groups is a sheaf on $inf(\Alg/A)^{op}$.
It was shown in [3,\S 5], that the groups $H_{inf}^*(A,\Cal G)$ as
defined above  are
the sheaf cohomology groups of the sheaf obtained by restriction of $\Cal
G$ to $inf(\Alg/A)$. In particular they are independent of the choice
of the presentation $A=R/I$. We also consider the hypercohomology of
sheaves of (cochain) complexes and of simplicial sets. If 
$\Cal G:\dots @>\delta>> \Cal G^n@>\delta>>\Cal G^{n+1}@>\delta>>\dots$ 
is a cochain complex
of sheaves, we write $\Bbb H_{inf}^*(A,\Cal G)$ for the cohomology
of the total complex:
$$
Tot_\pi^n=\prod_{p+q=n}C^p(R,I,\Cal G^{q}) 
$$
of the bicomplex $(p,q)\mapsto C^p(R,I,\Cal G^q)$. For a sheaf $X$ of
fibrant simplicial sets,
we put 
$\Bbb H_{inf}^r(A,X)=\pi_{-r}\holi_{\Delta\times\Bbb N^{op}}((m,n)
\mapsto XCyl^m(R,I)_n)$; the definition can be extended to Kan spectra as in 
[14]. 
\bigskip 
\remark{Other variants of infinitesimal cohomology 1.4} For a
unital algebra $A$, we can consider, in addition to the infinitesimal site
defined above, the indiscrete site on $inf(\Alg_1/A)^{op}$, where
$inf(\Alg_1/A)\subset inf(\Alg/A)$ is the subcategory of unit preserving
maps of unital algebras. However if $\Cal G:inf(\Alg/A)@>>>Ab$ is a
functor, then the sheaf cohomology groups $H^*(inf(\Alg_1/A),\Cal G)$
agree with those defined above. Indeed, both $H^*(inf(Alg_1/A),\Cal G)$
and $H^*(inf(Alg/A),\Cal G)$ can be computed as the cohomology of
$C(R,I,\Cal G)$ where $A=R/I$ is any presentation of $A$ as a
quotient of a unital quasi-free algebra. Similarly if $A$ is unital and
commutative, the cohomology of a functor $\Cal G:inf(\Comm/A)@>>>Ab$ is
the same whether considered as a sheaf on $inf(\Comm/A)^{op}$ or
$inf(\Comm_1/A)^{op}$; they are both computed as the cohomology of the
cosimplicial abelian group $[n]\mapsto \limo\Cal G(S^{\otimes
n+1}/\ker(S^{\otimes n+1}\fib A)^m)$ where $S\fib A$ is any surjective map
from a smooth unital algebra $S$, and $S^{\otimes n}\fib A$ is the induced
map. On the other hand, it follows from the discussion on page 337 of
[8] that, for $A$ commutative and unital, the cohomology of the additive
group functor $inf(\Comm_1/A)\owns (B\fib A) \mapsto B\in Ab$ agrees with
the cohomology of the structure sheaf on Grothendieck's infinitesimal site
$Spec(A)_{inf}$. Thus all three commutative infinitesimal cohomologies of
a commutative unital algebra agree. However these do not agree with the
noncommutative infinitesimal cohomology of the same algebra. In fact, we
shall show further below that for any algebra $A$, the cohomology of the
sheaf $O/[O,O]:inf(Alg/A)@>>>Ab$, $(B\fib A)\mapsto B/[B,B]$ is periodic
cyclic homology $HP_*(A)$. Since on the other hand, for a finite type
algebra over a field we have $HP_*(A)=\prod_{n=0}^\infty H^{2n+*}
(SpecA_{inf},O)$ (cf. [9]), we can see that the commutative and
noncommutative infinitesimal cohomologies do not agree in general.
Herefrom, we shall always use $H_{inf}$ to refer to the cohomology on the
noncommutative site. 
\endremark 
\bigskip
\subhead{2. Flat connections}\endsubhead 
\bigskip
\subsubhead{Infinitesimal modules, crystals and stratifications 2.0}
\endsubsubhead Let $A$ be an algebra. By
a right, {\it infinitesimal module} over $A$ or
$\tilde{O}_{inf(\Alg/A)^{op}}$-module we mean a functor
$M:inf(\Alg/A)@>>>Ab$ such that $M_p\in\Mod-\tilde{B}$ for each
$inf(\Alg/A)\owns p:B\fib A$. We say that the infinitesimal module $M$ is
a {\it crystal} on the infinitesimal site if for each map $\alpha:(p:B\fib
A)@>>>(p':B'\fib A)\in inf(\Alg/A)$ the $\tilde{B'}$-homomorphism
$\rho_\alpha:\alpha(M_p):=M_p\otimes_{\tilde{\alpha}}\tilde{B'}@>>>M_{p'}$
induced by $M(\alpha)$ is an
isomorphism. A {\it stratification} on a right $\tilde{A}$-module $E$ is a
sequence of isomorphisms $\theta_n:\faz^0_n(E)\iso
\faz^1_n(E)\in\Mod-Cyl^1(A)_n$ ($n\ge 0$); here $\faz^i_n$ is the $i$-th
coface map in the cosimplicial algebra $Cyl^*(A)_n$. The $\theta$'s are
subject to the following prescriptions: we must have $\theta_1=1$, the
identity of $E$; the sequence $\{\theta_n\}$ must be compatible with the
projection maps $Cyl^1(A)_{n+1}@>>>Cyl^1(A)_n$, and finally, each
$\theta_n$ must satisfy the cocycle condition: $$
\faz^2_n(\theta_n)\faz^0_n(\theta_n)=\faz^1_n(\theta_n)\tag{14} $$ For
example if $M$ is a crystal then the sequence
$\theta_n(M)=\rho_{\faz^1_n}^{-1}\rho_{\faz^0_n}^{}$ is a stratification
on $E(M)=M_{A@>id>>A}$. For a quasi-free algebra $R$, the 
assignment $M\mapsto (E(M),\theta(M))$, is an equivalence of categories
between the category of crystals of $\tilde{R}$-modules and natural
transformations and the category of pairs $(E,\theta)$ of an
$\tilde{R}$-module together with a stratification $\theta$, where a map
$(E,\theta)@>>>(E',\theta')$ is a module homomorphism $f:E@>>>E'$ such
that $\theta_n\faz^0_n(f)=\faz^1_n(f)\theta_n$ for all $n\ge 1$. An
inverse for this functor was constructed in [3,7.2]; here we write 
$E_{crys}$ for this
inverse equivalence. If $M$ is a crystal, then the pro-cosimplicial group
$MCyl^*(R)$ of 1.3 above can be described in terms of $E$ and $\theta$
as follows. We have $MCyl^0(R)=E$ and $MCyl^m(R)=\faz^m\dots\faz^1(E)$ for
$m\ge 1$.  Codegeneracies are the maps
$id_E\otimes\mu^i:\faz^m\dots\faz^1(E)@>>>\faz^{m-1}\dots\faz^1(E)$ ($0\le
i\le m$). Similarly, for $1\le j\le m+1$, $id_E\otimes\faz^j$ is the
$j$-th coface map $\faz^m\dots\faz^1(E)@>>>\faz^{m+1}\dots\faz^1(E)$.  
The $0$-th coface map $E@>>>\faz^1(E)$ is the composite
$\faz^0_\theta:E@>>>\faz^0(E)@>\theta>>\faz^1(E)$. For $m\ge 1$
and $e\otimes\omega\in\partial^m_n\dots\partial^1_n(E)$ we have
$\partial^0_\theta(e\otimes\omega)=
\partial^0_\theta(e)\otimes\partial^0(\omega)$.
Upon
completion we obtain a complex $C(E,\theta)$ whose cohomology is the sheaf
cohomology of the crystal $M$ with which we started.

We have shown how one can go from crystals on a quasi-free algebra 
to stratifications and interpret
the sheaf cohomology of the former in terms of the latter. Below we shall
show how one can go from stratifications to flat connections. First we
recall the definition of curvature given in [11].
Given a connection $\nabla:E@>>>E\otimes_{\tilde{A}}\Omega^1(A)$ 
on a $\tilde{A}$-module $E$, extend it to a map 
$\nabla:E\otimes_{\tilde{A}}\Omega^*(A)@>>>E\otimes_{\tilde{A}}
\Omega^{*+1}(A)$
by the rule $\nabla(e\otimes\omega)=\nabla(e)\omega+ed\omega$. The 
{\it curvature}
of $\nabla$
is the map $\nabla^2:E\otimes_{\tilde{A}}\Omega^*(A)@>>>E\otimes_{\tilde{A}}
\Omega^{*+2}(A)$. We say that $\nabla$ is {\it flat} if $\nabla^2=0$. One
sees immediately that flatness is equivalent to the vanishing of the
degree $2$ component $\nabla^2:E@>>>E\otimes_{\tilde{A}}\Omega^2(A)$. 
The pairs $(E,\nabla)$ of a right $\tilde{A}$ module and a flat connection 
form a category, where a map $(E,\nabla)@>>>(F,D)$ is a homomorphism
 of modules $\alpha:E@>>>F$ with 
$D\alpha=\alpha\otimes id_{\Omega^1(A)}\nabla$. 
\bigskip
\proclaim{Proposition 2.1} Let $R$ be a quasi-free algebra. With the 
definitions of 2.0 above, let          
$(E,\theta)$ be a stratified $\tilde{R}$-module and 
$\faz^0_{\theta,1}:E@>>>\faz^1_1(E)$ the coface map.
Then the composite
$\nabla:E@>\faz^0_{\theta,1}>>\faz^1_1(E)=E\oplus 
E\otimes_{\tilde{A}}\Omega^1(A)\fib E\otimes_{\tilde{A}}\Omega^1(A)$ 
is a flat connection. Moreover the assignment $(E,\theta)\mapsto (E,\nabla)$,
$f\mapsto f$ is an isomorphism between the categories of stratified
modules and of flat connections.
\endproclaim

\demo{Proof} That $\nabla$ is a connection is immediate, as is that a 
homomorphism of stratifications is a morphism of modules with a connection. 
The proof that $\nabla$
is flat will be easy once we have established the identity \thetag{15} below;
to
prove the latter, proceed as follows. Consider the map 
$\faz^0_\theta:E@>>>E\otimes_{\partial_1} Cyl^1(R)\cong
E\otimes_{\tilde{R}}\Omega(R)/\Omega^+(R)^\infty$; write 
$\faz^0_\theta=\sum_{n=0}^\infty D_n$ where 
$D_n(E)\subset E\otimes_{\tilde{R}}\Omega^n(R)$.
Thus $D_0=1$ and $D_1=\nabla$. Similarly, the map 
$\faz^0_\theta:\faz^1(E)@>>>\faz^2\faz^1(E)$, decomposes as
$\faz^0_\theta=\sum_{n=0}^\infty\sum_{i=0}^nD_i\otimes P_{n-i}$,
where the $P$'s are the components of the map 
$\faz^0:Cyl^1(R)@>>>Cyl^2(R)$. The identity 
$\faz^0_\theta\faz^0_\theta=\faz^1\faz^0_\theta$ is
equivalent to the sequence of identities:
$$
D_n(e)\otimes v_2^n-D_n(e)\otimes v_1^n-D_n(e)\otimes (v_2-v_1)^n=
\sum (D_i\otimes P_j)
(D_k(e)\otimes v_1^k)\tag{15}
$$
where $e\in E$, $n\ge 1$, and the sum is taken over all $i\ge 0$ and $j,k\ge 1$
such that $i+j+k=n$. The case $n=1$ of the identity \thetag{15} is 
tautological. That $\nabla$ is flat follows straightforwardly from the case
$n=2$ of the identity above, using the formula \thetag{9}. Thus
$(E,\theta)@>>>(E,\nabla)$ is a functor from stratified modules to modules
with a flat connection. Next we construct a functor going in the
opposite 
direction.
Consider the map $t:\Omega^1(R)@>>>Cyl^1(R)$, $adx\mapsto a(\faz^0(x)-x)$.
Define a map $\faz^0_\nabla:E@>>>\faz^1(E)$, $\faz^0_\nabla(e)=e+(id\otimes t)\nabla(e)$,
and put $\theta:\faz^0(E)@>>>\faz^1(E)$, 
$\theta(e\otimes\omega)=\faz^0_\nabla(e)\faz^0(\omega)$. We must see now that
$\theta$ is a
stratification.
One checks that $\faz^0_\nabla(ea)=\faz^0_\nabla(e)\faz^0(a)$ ($e\in E, a\in R$),
from which it follows that $\theta$
is a well-defined homomorphism. Of the conditions for a stratification
only the cocycle condition is not immediate. To prove the latter, one does
as follows. One checks that \thetag{14} is equivalent to the 
cosimplicial identity $\faz^0_\theta\faz^0_\theta=\faz^1\faz^0_\theta$.
In turn, since $\faz^0_\theta=\faz^0_\nabla:E@>>>\faz^1(E)$, the latter
condition is equivalent to the sequence of identities \thetag{15} with
as $D_n$ the $n-th$ component of $\faz^0_\nabla$. A straightforward induction
shows that $\thetag{15}$ holds. On the other hand it is clear from the
definition of $D_n$ that if $f:(E,\nabla)@>>>(E',\nabla')$ is a homomorphism
of modules with a connection, and $\theta$ and $\theta'$ are the 
stratifications
constructed from them, then $\partial^1(f)\faz^0_{\theta}=\faz^0_{\theta'}f$. It follows 
that $f$ is a map of stratified modules 
$(E,\theta)@>>>(E',\theta')$. Thus we have a functor from flat connections
to stratifications. It is clear than if we start off with a flat connection,
then take the stratification constructed from it and then go back by
the functor of the proposition, we end up with the same connection we
started with. That the other composition is the identity is a consequence
of the recursive nature of the formula \thetag{15}; indeed, the latter
formula says that all of the $D_n$ --and therefore the whole stratification--
are determined by $D_1=\nabla$.\qed
\enddemo
\bigskip
\proclaim{Corollary 2.2} For a quasi-free algebra $R$, the
categories of $\tilde{R}$ modules with a flat
connection and of infinitesimal crystals are equivalent.
\endproclaim
\demo{Proof} It follows from the proposition and the equivalence between
the categories of stratified modules and crystals discussed in 2.0 
above.\qed
\enddemo
\bigskip
\proclaim{Theorem 2.3}Let $R$ be a quasi-free algebra, $E$ a right 
$\tilde{R}$-module, $\nabla:E@>>>E\otimes_{\tilde{R}}\Omega^1(R)$ a
flat connection in the sense of 2.0, and $E_{crys}$ the infinitesimal
crystal corresponding to $(E,\nabla)$ under the category equivalence
of corollary 2.2. Then there is an isomorphism:
$$
H^*_{inf}(R,E_{crys})\cong H^*(E\otimes_{\tilde{R}}\Omega(R),\nabla)
$$
\endproclaim
\bigskip
The proof of the theorem above uses the following variant of the perturbation
lemma of [12]:
\bigskip
\proclaim{Lemma 2.4} Let $C=(\{C^n\}_{n\ge 0},\delta)$ be a nonnegatively 
graded
cochain complex. Assume there is a decomposition $C^n=\prod_{p=0}^\infty C^n_p$
($n\ge 0$) and that $\delta$ preserves the associated filtration, i.e. 
$\delta(\prod_{p\ge q}^\infty(C_p^*))\subset\prod_{p\ge q}^\infty(C_p^*)$
for
all $q\ge 0$. Write $\delta=\sum_{n=0}^\infty\delta_n$ where 
$\delta_n(C^*_*)\subset (C^{*+1}_{*+n})$. Suppose $h_0:C^*_*@>>>C^{*-1}_*$ is
a contracting homotopy for $\delta_0$ in the sense that 
$1=h_0\delta_0+\delta_0h_0$. Define recursively:
$$
h_n=-\sum_{i=0}^{n-1}h_0\delta_{n-i}h_i
$$
Then $h=\sum_{n=0}^\infty h_n$ is a contracting homotopy for $\delta$.
\endproclaim
\demo{Proof} Straightforward induction.\qed\enddemo
\bigskip
\demo{Proof of Theorem 2.3} Consider the complex $C=C(R,0,E_{crys})$ of
\thetag{13}; this
is the cochain complex associated to the completion of the pro-cosimplicial 
abelian group $C(E,\theta)$ of  
2.0. We shall define a map going from the normalized complex $\bar{C}$ 
to the complex $(E\otimes_{\tilde{R}}\Omega(R),\nabla)$ and 
show it is a quasi-isomorphism. First we need a better 
description of the complex $\bar{C}$. We have a decomposition
$C^m=\prod_{n=0}^\infty C^m_n$, where 
$C^m_n:=E\otimes_{\tilde{R}}\Omega^n(R)\otimes T^n(V^m)$; cofaces and
codegeneracies are as described in 2.0 above. Thus all cofaces but the
$0$-th
coface act only on the $T(V)$ part, i.e. are of the form $id\otimes \faz^i$
where $\faz^i$ is the coface of the completed tensor product cosimplicial 
vectorspace $T(V):m\mapsto \prod_{n=0}^\infty T^n(V^m)$. 
The $0$-th coface map decomposes as an infinite sum of homogeneous
components of nonnegative degrees of which the $0$-th degree component is 
of the form $id\otimes \faz_0$ where $\faz_0$ is the coface of $T(V)$.
Hence the decomposition above is preserved by normalizing and we have 
$\bar{C}^m=\prod_{n=0}^\infty\bar{C}^m_n$
 where $\bar{C}^m_n=E\otimes_{\tilde{R}}\Omega^n(R)\otimes \bar{T}^n(V^m)$.
Here $\bar{T}^n(V^m)$ is the normalization of the cosimplicial tensor product
of the cosimplicial vectorspace of 1.1. 
By definition, $\bar{T}^n(V^m)$ is the
quotient of $T^n(V)$ by the subspace generated by the sum of the images of
all cofaces but the $0$-th coface; one checks this subspace is that generated
by all those words of length
$n$ on the basis elements $v_1,\dots, v_m$ of $V^m$ in which at least one
of the $v_i$ is missing. Hence each $\bar{T}^n(V^m)$ is isomorphic to the free
vectorspace on the set of all surjective maps going from a set of $n$
elements to a set of $m$ elements. In particular $\bar{T}^n(V^m)=0$ for
$m>n$, and $\bar{T}^n(V^n)\cong k[S_n]$, the group algebra of the symmetric
group. The coboundary in $\bar{C}$ is the map $\delta$ induced by $\faz^0$;
it admits a decomposition as that of the lemma above, where $\delta_0$ is
of the form $id\otimes\partial'_0$, with $\partial'_0$ the map induced
by the $0$-th coface of $T(V)$ on $\bar{T}(V)$. Consider the map
$p:k[S_n]\fib k$,
$\sigma\mapsto sg(\sigma)$
where $sg$ is the permutation sign. One checks that 
$\ker p=\partial'_0(T^n(V^{n-1}))$; thus if we regard $p$ as a map of
cochain
complexes $\bar{T}^n\fib k[-n]$, we see it induces an 
isomorphism at the $H^n$ level. I claim it is in fact a quasi-isomorphism,
i.e. we have $H^q(\bar{T}^n(V))=0$ for $q<n$. The claim follows from the
fact that
$\bar{V}=k[-1]$, the K\"unneth formula for the cohomology of cochain complexes
and the fact that --as can be seen through the combined application of the 
Eilenberg-Zilber theorem and the Dold-Kan normalization theorem-- the 
normalized chain complex of the tensor product $n\mapsto A^n\otimes B^n$
of two cosimplicial vectorspaces $A$ and $B$ is quasi-isomorphic to the 
tensor product of the normalized cochain complexes. Thus the claim is proved,
and we may regard $p:\bar{T}^n(V)\fib k$ as a free resolution of $k$;
in particular the cochain complex $\ker p^*$ is contractible. One checks that,
upon
tensoring over $k$ with $E\otimes_{\tilde{R}}\Omega(R)$
 and completing, we obtain a cochain map
$1\otimes p:(\bar{C},\delta)@>>>(E\otimes_{\tilde{R}}\Omega(R),\nabla)$. 
By what we have just seen, the kernel of this map satisfies the hypothesis
of the lemma above. This proves the theorem.\qed
\enddemo
\bigskip
\subhead{3.Periodic homology v. infinitesimal cohomology}\endsubhead
The purpose of this section is to prove the following:
\bigskip
\proclaim{Theorem 3.0} Let $k$ be a field of characteristic zero, 
$A$ a $k$-algebra. Then there is an isomorphism:
$$
HP_*(A)\cong H^*_{inf}(A,O/[O,O])
$$
\endproclaim
\demo{Proof} Write $A=R/I$, with $R$ quasi-free. The proof has two parts. 
First we prove that 
$H^*_{inf}(A,O/[O,O])$ equals the cohomology of the complex:
$$
(\lim_n {(\frac{\Omega^*}{\Cal{F}^n+b\Omega^{*+1}})}_\kappa,d)\tag{16}
$$
where $\Omega^*=\Omega^*(R)$, 
$\Cal F$ is the adic filtration of the ideal $I$, as defined in 1.2 above,
and $\kappa$ is the Karoubi operator. The map $\kappa$ verifies 
$\kappa^m\equiv 1$ on $\Omega^m$ $\mod$ $(\Cal{F}^n+b\Omega^{m+1})$;
the 
subscript
in \thetag{16} indicates coinvariants with respect to this action of 
$\Bbb Z/m$.
The second part consists of proving that the latter complex computes 
$HP_*(A)$. Now to the first part. Consider the 
inverse systems of normalized cochain complexes $\{\bar{C}_n:n\ge 0\}$
and $\{\bar{D}_n:n\ge 0\}$ 
associated with the 
pro-cosimplicial groups $\{Cyl(R,I)_n/[R,Cyl(R,I)_n]: n\ge 0\}$ and
 $\{Cyl(R,I)_n/[Cyl(R,I)_n,Cyl(R,I)_n]: n\ge 0\}$. Thus, by \thetag{11},
$\bar{C}_n^m=\oplus_{i=m}^n\bar{C}_{n,i}^m$, where, for $\bar{T}$ as in the 
proof of Theorem 2.3,
$\bar{C}_{n,i}^m=\frac{\Omega^i}{F^{n-i}\Omega^p+b\Omega^{i+1}}
\otimes\bar{T}^i(V^m)$. We write 
$\delta_n=\sum_{i=0}^n\delta_{n,i}$ for the coboundary of $\bar{C}_n$.
 A similar argument as that in the proof of assertion
c) on page 30 of [11] shows that $\bar{D}_n^m=\oplus_{i=m}^n\bar{D}_{n,i}^m$
where $\bar{D}_{n,i}^m=(\bar{C}_{n,i}^m)_{\Bbb Z/i}$ is the coinvariants
with respect
to the action given on a generator $\zeta$ of $\Bbb Z/i$ by 
$\zeta(\omega\otimes x)=\kappa'(w)\otimes \zeta(x)$. Here 
$\kappa'(a_0da_1\dots da_i)=da_ia_0da_1\dots da_{i-1}$ is Karoubi's
operator without sign, and $\zeta(x_1\dots x_i)=x_ix_1\dots x_{i-1}$.
Proceeding as in the proof of theorem
2.3 we see that there is an inverse system of cochain maps 
$\pi_n:\bar{C}_n\fib \Omega^*/\Cal F^n$ which are $0$ in codimension
$>n$ and are given in codimension $\le n$ by the composite
of the projection 
$\bar{C}^m_n@>>>\bar{C}^m_{n,m}=\frac{\Omega^m}{F^{n-m}\Omega^m}\otimes k[S_m]$
and the map $(1\otimes p)$, where $p$ is as in the proof of theorem 2.3.
Thus the restriction $p_i$ of $p$ to each of the complexes $\bar{T}^i(V)$ is
a quasi-isomorphism of $k[\Bbb Z/i]$-modules. Since $char(k)=0$, we
therefore
have that $\ker p_i$ is contractible as a complex of $k[\Bbb 
Z/i]$-modules,
i.e. there is a contracting homotopy which commutes with the action
of $\Bbb Z/i$. This contracting homotopy is preseved upon tensoring with 
$\frac{\Omega^i}{F^{n-i}\Omega^{i}+b\Omega^{i+1}}$ and gives a map of inverse
systems 
$\ker(\pi_n)^*@>>>\ker(\pi_n)^{*-1}$ which commutes with the action of the
cyclic groups and is a contracting homotopy for $\delta_{n,0}$. Now apply
Lemma 2.4 to obtain an inverse system of action preserving contracting
homopies
for $\ker\pi_n$. Taking coinvariants and then inverse limit, we obtain
a contracting homotopy for $\ker\limu\pi'_n$, where $\pi'_n$ is the map
induced by $\pi_n$ at the coinvariant level. To finish
the first part of the proof it suffices to show that $\pi'=\limu\pi'_n$
is surjective. To see this, consider the map 
$s_n:\Omega^m/F^{n-m}\Omega^m @>>> C^m_n$, 
$\omega\mapsto\omega\otimes(1/m!\sum_{\sigma\in S_m}\sigma)$. We see that
$s_n$ is a right inverse for $\pi_n$, that it commutes with the action
of the cyclic groups --although not with coboundaries-- and that the 
family $\{s_n\}$ is a map of inverse systems. Hence if $s'_n$ is the map
induced by $s_n$ at the coinvariant level and $s'=\limu s'_n$, then
$\pi's'=1$ and $\pi'$ is surjective. 

Next we turn to the second part of the proof. First of all, note that,
because
the map $1-\kappa$ is invertible on the second summand of the harmonic
decomposition $\Omega=P\Omega\oplus P^\perp\Omega$ of 
[6], we have an isomorphism of cochain complexes
$(P\frac{\Omega}{\Cal F^n}_\natural,d)\cong
((\frac{\Omega^m}{F_{n-m}\Omega^m+b\Omega^{m+1}})_\kappa,d)$ for each $n\ge 0$.
Here
$P\frac{\Omega^m}{\Cal 
F^n}_\natural=\frac{P\Omega^m}{P(b\Omega^{m+1}+F_{n-m}\Omega^m)}$;
we abuse notation and write $d$ for the map it induces at the 
$_\natural$ level.
We want to compare the latter cochain complex with the supercomplex $X$
of [5], which is known to compute the periodic cyclic groups.
Since these complexes live in different categories we introduce $\tau X$,
a nonnegatively 
graded cochain complex version of the supercomplex $X$. We put 
$\tau X^0=\Omega^0_\natural$, and for $n\ge 1$, 
$\tau X^{2n}=\Omega^0=X^{even}$, $\tau X^{2n-1}=\Omega^1_\natural=X^{odd}$.
Abusing notation and writing $\Cal F^\infty$ for the image of the filtration
\thetag{12}, we have $H^*(\tau X/\Cal F^\infty)=H^*(X/\Cal F^\infty)$.
It follows from [2, 2.3] that the filtration $\Cal F$ is equivalent to
the
filtration induced by $K^n=\ker(\Omega\fib \Omega(R/I^n))$;
in turn it was shown in [5] that for $R=TA$ and $I=JA$, the latter 
filtration is equivalent to the Cuntz-Quillen filtration. Thus 
$H^*(X/\Cal F^\infty)=HP_*(A)$. We shall presently define inverse
homotopy equivalences of cochain complexes 
$P\Omega_\natural\rightleftarrows\tau X$.
First we introduce the following two operators; let 
$h:\Omega^*@>>>\Omega^{*-1}$, $h(\omega da)=(-1)^{|\omega|}\omega a$, 
and $\nabla:\Omega^*@>>>\Omega^{*+1}$, 
$\nabla(a_0da_1\dots a_n)=-a_0\phi(a_1)da_2\dots da_n$ where $\phi$ is as
in 1.1 above. One checks the following properties of $h$ and $\nabla$. We have:
$$
(\nabla b+b\nabla)\omega=\omega, \quad\nabla(uv)=\nabla(u)v+(-1)^{|u|}udv,
\quad\nabla\Cal F^n\subset\Cal F^n\tag{17}
$$ 
for $\omega\in\Omega^{2+*}$, $u,v\in \Omega^{1+*}$ and also:
$$
h^2=0\quad hd+dh=1\quad h\Cal F^n\subset\Cal F^{n-1}\tag{18}
$$
Below we define maps $\alpha:P\Omega_\natural@>>>\tau X$,
$\beta:\tau X@>>>P\Omega_\natural$, $\gamma:\tau X@>>>\tau X[-1]$,
and $\epsilon:P\Omega_\natural@>>>P\Omega_\natural[-1]$. One checks,
using \thetag{17}, \thetag{18},\thetag{12} and induction, that all these maps
are continuous in the adic topology,
that $\alpha$ and $\beta$ are cochain maps and that $\gamma$ and $\epsilon$
are homotopies $\alpha\beta@>>>1$ and $\beta\alpha@>>>1$. For $j=0,1$,
$1\le i\le q$, and $n=2q+j$, we write
$c_{n,i}=n(n-2)(n-4)\dots (2i+j)$, and 
$c_{n}=c_{n,1-j}$.  
Let $\alpha^0$, $\alpha^1$, $\beta^0$ and $\beta^1$ be identity maps, 
$$
\alpha^n=c_n(Phb)^{q},\quad \beta^{2q+1}=\frac{1}{(2q+1)!}
(P(\nabla B+B\nabla))^{q-1}P,\quad
\beta^{2q}=(2q+1)\beta^{2q+1}\nabla B$$
Here and below the superscrypts on expressions between parenthesis denote
powers; thus for instance $(Phb)^{q}=Phb\circ\dots\circ Phb$ ($q$ times).
Put also $\gamma_i=0$ for $i\le 1$,
$\epsilon_l=0$ for $l\le 2$, and: 
$$
\gamma^{2n+1}=hP(b\nabla -1)-\sum_{i=1}^n
\alpha^{2i}(PhP\nabla b)\beta^{2i+1},
\qquad\gamma^{2n}=-\sum^n_{i=1}\alpha^{2i-1}PhP\nabla b\beta^{2i}$$
$$
\epsilon^{2q+1}=(2q+1)\sum_{i=1}^q\frac{1}{c_{2q,i}}
(P(\nabla B+B\nabla))^{q-i}P\nabla(Phb)^{q-i+1}$$
$$\epsilon^{2q}=\sum_{i=2}^q\frac{1}{c_{2q-1,i-1}}
(P(\nabla B+B\nabla))^{q-i}P\nabla(Phb)^{n-i+1}.\qed$$
\enddemo
\bigskip
\remark{Remark 3.1} The cohomology isomorphism 
$H^n(P\frac{\Omega}{\Cal F^\infty}_\natural)
\cong HP_{n}(A)\cong H^{n+2}(P\frac{\Omega}{\Cal F^\infty}_\natural)$ ($n\ge 0$) is
induced by
the cochain map $s':=PhbN:P\Omega_\natural@>>>P\Omega_\natural[-2]$, 
where $N(\omega)=|\omega|\omega$. This follows from the fact that, for the homotopy
equivalence of the proof above, we have 
$\alpha^{n}{s'}^{n+2}=\alpha^{n+2}$.
\endremark
\bigskip
\subhead{4. A  new presentation of the Chern character}\endsubhead
The title of this section refers to the Connes-Karoubi map:
$$
K_0(A)@>ch_0>>HP_0(A)\tag{19}
$$
We shall demonstrate here how this character fits into the infinitesimal
cohomology framework. In 4.1 below we shall construct a character 
$ch^{inf}:K_0(A)@>>>H^2_{inf}(A,K_1)$, where $K_1$ is Bass's. This 
character is defined independently of the characteristic  of 
$k$; we show in 4.4 that if $char(k)=0$, then our character is the same 
as \thetag{19}, up to a factor of $2$. The key to proving this agreement 
result is provided by the following:
\bigskip
\proclaim{Proposition 4.0} Let $A$ be an algebra over a field $k$ of
characteristic zero; then there is a natural
isomorphism:
$$\nu:H^n_{inf}(A,K_1)\cong HP_n(A)\quad(n\ge 2)$$
{\smc Note:} The groups $H^n_{inf}(A,K_1)$ for $n=0,1$ are discussed in 
section 5.
\endproclaim
\demo{Proof} Write $R/I=A$ with $R$ quasi-free. By virtue of \thetag{11}
we have, for each $n\ge 1$, an exact sequence of cosimplicial algebras:
$$
0@>>>\oplus_{i=1}^{n-1}\frac{\Omega^i(R)}{F^{n-i}\Omega^i(R)}\otimes T^i(V^*)@>>>
Cyl^*(R,I)_n@>\pi_n>>\frac{R}{I^{n}}@>>>0
$$
here $\frac{R}{I^n}$ is the constant cosimpicial algebra. We note that we
also have a natural inclusion map $\frac{R}{I^n}\subset Cyl^*(R,I)_n$ which 
is an algebra
homomorphism (although not a cosimplicial map), and also a right inverse
for $\pi_n$ and a map of inverse systems. 
Thus if $\Cal G$ is any functor going from algebras
to abelian groups, the sequence 
$$
0@>>>\lim_n\ker \Cal G(\pi^*_n)@>>>\lim_n\Cal G(Cyl^*(R,I)_n)@>>>
\lim_n \Cal G(\frac{R}{I^n})@>>>0\tag{20}
$$
is exact. The proposition follows from the isomorphism 
$\ker K_1(\pi^*_n)\cong\ker HC_0(\pi^*_n)$ of [13], theorem 3.0 above, 
and the long cohomology sequences
associated to \thetag{20} for $\Cal G=K_1$, and for
$\Cal G=HC_0=O/[O,O]$.\qed
\enddemo
\bigskip
\definition{The infinitesimal Chern character 4.1}
The construction has two
parts; first we define, for any algebra $A$, a character 
$ch^{strat}:K_0(A)@>>>H^2_{strat}(A,K_1)$ taking values in the cohomology
of the stratifying
site of [3,5.6], i.e. in the cohomology of the complex 
$\limu K_1(Cyl(A)_n)$.
Note that for $A$ quasi-free this suffices, as 
$H^*_{strat}(A,-)=H^*_{inf}(A,-)$ in this case. The second part consists of
lifting the construction above to the infinitesimal topology of an arbitrary
algebra. Now, to the first part. 
Suppose a finitely generated projective right $\tilde{A}$-module
$E$ is given; then I claim that there exists, 
for each $n\ge 1$,
an isomorphism of $\widetilde{Cyl}(A)_n$-modules
$\theta_n:\partial^0_n(E)\iso\partial^1_n(E)$ with $\theta_1=1$ which
is compatible with the projection maps $\partial^i_{n+1}(E)@>>>\faz^i_{n}(E)$
($i=0,1$).
This is best seen in terms of idempotent matrices. 
Let $p\in M_n(\tilde{A})$ be an idempotent matrix with $p\tilde{A}^n=E$;
we have 
to show that there exists a compatible family of invertible
matrices $u_n\in Gl_n(\widetilde{Cyl}(A)_n)$ with $u_1=1$ and
such that $u_n\faz^0_n(p)=\faz^1(p)u_n$. Since an idempotent
matrix is the same thing as a homomorphism $k@>>>M_n(\tilde{A})$, it
suffices to
show
the above for $A=k$ and $n=1$. In this particular case, a formula
for $u$ is given in Lemma 4.2 below.
The claim is thus proved. Next we consider
the map $\faz_n(\theta_n):=\faz^2_n(\theta_n)\faz^0_n(\theta_n)
\faz^1_n(\theta_n)^{-1}$; we observe that this  
is an automorphism of the $\widetilde{Cyl}^2(A)_n$-module
$\partial_2\partial_1(E)_n$,
and therefore has a class in $K_1(\widetilde{Cyl}^2(A)_n)$. It is immediate
from a well-known form of Whitehead's lemma --which we have included
below as lemma 4.3-- that the class we have constructed
is a $2$-cocycle of the complex $K_1(\widetilde{Cyl}^*(A)_n)$. 
Because the $\theta_n$ are
compatible, so is the sequence $\faz(\theta)=\{\faz_n(\theta_n)\}$ whence
it is a $2$-cocycle of the completed cochain complex 
$\limu K_1(\widetilde{Cyl}^*(A)_n)$. We write $ch^{strat}(E)$ for its 
cohomology class in $H_{strat}^2(A,K_1\circ\widetilde{\quad})=
H_{strat}^2(A,K_1(-)\oplus K_1(k))=H_{strat}^2(A,K_1)$. 
It follows by lemma 4.3 that if
$\theta'$ is another inverse system of isomorphisms as above, then the
automorphisms  
$\faz_n(\theta_n)\faz_n({\theta'}_n)^{-1}$ of $\faz^2_n\faz^1_n(E)$ and
$\faz^2_n(\theta_n{\theta'}^{-1}_n)\oplus\faz^1_n(\theta_n\theta^{-1}_n)^{-1}
\oplus\faz^0_n(\theta_n{\theta'}^{-1}_n)$ of 
$\partial^2_n\partial^1_n(E)^2\oplus\partial^0_n(E)$ have the same class
in $K_1(\widetilde{Cyl}^1(A)_n)$. Hence the definition of
$ch^{strat}(E)$ is independent of the choice of $\theta$. Next, note
that if
$f:E\iso F$ and $\alpha:\faz^0(E)@>>>\faz^1(E)$ are isomorphisms, then so
is
$\beta=\partial^1(f)\alpha\partial^0(f)^{-1}$; one checks, using lemma 4.3,
that $\faz_n(\beta)\equiv\faz_n(\alpha)$ in $K_1(\widetilde{Cyl}^1(A))_n$. 
Hence
$ch^{strat}$ is well-defined on isomorphism classes of finitely generated
projective $\tilde{A}$ modules. Furthermore it is a monoid 
homomorphism, because 
if $\alpha:\faz^0(E)@>>>\faz^1(E)$ and $\beta:\faz^0(F)@>>>\faz^1(F)$
are the isomorphisms used to define $ch^{strat}(E)$ and $ch^{strat}(F)$, 
then we can use $\alpha\oplus\beta$ for $E\oplus F$, and 
$\partial(\alpha\oplus\beta)=\partial(\alpha)\oplus\partial(\beta)$. 
Hence we have
a well-defined group homomorphism
$K_0(\tilde{A})@>ch^{strat}>>H^2_{strat}(A,K_1)$, and by composition with
$K_0({A})@>>>K_0(\tilde{A})$ or --since, as one sees
immediately,  $ch^{strat}(\tilde{A})=0$-- by passage to the quotient 
$\mod$ $K_0(k)$, also a homomorphism 
$K_0(A)@>ch^{strat}>>H^2_{strat}(A,K_1)$. 
This homomorphism is natural, 
as follows from the fact that lemma 4.2
gives a formula for an isomorphism $\theta:\faz^0(E)@>>>\faz^1(E)$,
whence also for $ch^{strat}(E)$.

Now we lift this construction to the infinitesimal site as follows.
Write $A=R/I$ where $R$ is quasi-free.
Given a projective right $\tilde{A}$-module lift it succesively
to an inverse system $\pi_{n}:E_{n}@>>>E_{n-1}$ where $E_1=E$, and
each $\pi_n$ is an $\tilde{R}/I^n$ homomorphism. That this can be done
is a consequence of the fact that $k$ is quasi-free; explcit formulas
for the succesive idempotents $e_n$ for each of the $E_n$ in terms of that
of $E$ are given in [5] for $char(k)=0$. Using as 
$\theta_n:\faz^0_n(E_n)@>>>\faz^1_n(E_n)$ the image through the homomorphism
induced by $e_n$ of the truncation of the element $u$ of Lemma 4.2,
we obtain an element 
$\{\partial_n(e_n(u_n))\}\in\limu Gl(\partial^2_n\partial^1_nE_n)$.
Taking the class of this element at each level, we obtain a $2$-cocycle
$\faz(e(u))$
of the complex
$\limu K_1(\widetilde{Cyl}^*(R/I^n)_n)\cong\limu K_1(\widetilde{Cyl}^*(R,I)_n)$.
The last isomorphism comes from the fact that, by [2, 2.3], the filtrations
$(q^*(R)+<I>)^\infty$ and $\{\ker(Q^n(R)\fib Q^n(R/I^n)_n\}$ are
equivalent.
We write $ch^{inf}(E)$ for the cohomology class of $\faz (e(u))$. After
what we have already seen, to show
that $E\mapsto ch^{inf}(E)$ is well-defined on isomorphism classes of 
projective modules and gives a group homomorphism 
$ch^{inf}:K_0(A)@>>>H^2_{inf}(A,K_1)$, it suffices to show that any
isomorphism of projective modules $\alpha_1:E@>>>F$ can be lifted to an inverse
system of isomorphisms $\alpha_n:E_n@>>>F_n$ of the chosen liftings
of $E$ and $F$. But this is immediate from Lemma 4.2 and the fact that
the algebra $k[t,t^{-1}]$ is quasi-free.
\enddefinition
\bigskip
\proclaim{Lemma 4.2} Let $k$ be the ground field; write $e$ for its unit
element. Then:
\smallskip
\item{1)} There is an isomorphism of inverse systems
 $Cyl^1(k)\iso\frac{\Omega(k)}{\Omega^+(k)^\infty}$ such that
the canonical coface maps $\partial^i:k@>>>Cyl^1(k)$ $i=0,1$ are
$\faz^1(e)=e$ and
$$
\faz^0(e)=e+de+\sum_{n=1}^\infty C_n(1-2e)de^{2n}\tag{21}
$$
Here $C_1=1$ and $C_n=\sum_{i=1}^{n-1}C_iC_{n-i}$ is the $n+1$-th
Catalan number (cf. [10]).
\smallskip
\item{2)} The following element 
$u\in\frac{\widetilde{\Omega}(k)}{\Omega^+(k)^\infty}$
verifies the identity $u\partial^0(e)=\partial^1(e)u$. 
$$
u=1+\sum_{n=1}^\infty C_n(2e-1)de^{2n-1}\tag{22}
$$

\endproclaim
\demo{Proof}
One checks that the $1$-cocycle $\phi:k@>>>\Omega^2(k)$, $\phi(e)=(1-2e)de^2$
satisfies $-\delta(f)=d\cup d$. Part 1) of the lemma follows from this and 
the recursive formula (5') of [2].
Write $D_n(e)$ for the term of degree $n$ of \thetag{21}; for
$u=1+\sum_{n=1}^\infty \omega_n$, the condition $u\faz^0(e)=\faz^1(e)u$
is equivalent to the sequence of identities:
$$
[e,\omega_n]=D_n(e)+\sum_{i=1}^{n-1}\omega_iD_{n-i}(e)\tag{23}
$$
One checks that for $\omega_{2n-1}=C_n(2e-1)de^{2n-1}$, $\omega_{2n}=0$,
the identity above is satisfied. \qed
\enddemo
\bigskip
The following well-known version of Whitehead's Lemma was used in 4.1
above.
\bigskip
\proclaim{Whitehead's Lemma 4.3} Let $E$ and $F$ be right
$\tilde{A}$-modules,
and let $\phi:E\iso F$ and $\psi:F\iso E$ be module isomorphisms. Then the
following automorphisms:
$$
\left(\matrix\psi\phi &0 \\
        0 &1_{F\oplus E\oplus F}\endmatrix\right)\sim 
\left(\matrix 1_{E\oplus F\oplus E} &0\\ 
              0 & \phi\psi\endmatrix\right)
$$
are in the same commutator class of the group $Aut(E\oplus F\oplus E\oplus F)$.
\endproclaim
\bigskip
\proclaim{Theorem 4.4} Let $k$ be a field of characteristic zero, 
$A$ an algebra over $k$ and $\nu$ and $ch^{inf}$ be the homomorphisms of
4.0 and 4.1 above. Then:
$$
\nu\circ ch^{inf}=2ch_0
$$
\endproclaim
\demo{Proof} We keep the notations of 4.1. We shall show that,
for $E=e_1\tilde{A}^s$, the sequence of isomorphisms:
$$
\alignat2
H^2_{inf}(A,K_1)&\cong H^2_{inf}(A,HC_0)\qquad\text{by 4.0}\\
&\cong H^2(\lim_n\frac{P\Omega^*}{P(\Cal F^n+b\Omega^{*+1})}) 
\qquad\text{by the proof of 3.0}\\
&\cong
H^0(\lim_n\frac{P\Omega^*}{P(\Cal F^n+b\Omega^{*+1})})\qquad\text{via 3.1}\\
&\cong H^0(\hat{X}(R,I))\qquad\text{by the proof of 3.0}\\
&\cong HP_0(A)\quad\qquad\text{by [5]}\\
\endalignat
$$
maps $ch^{inf}(E)$ to twice the class of $\{e_n\}$ in $H^0(\hat{X}(R,I))$; by 
[5, 12.4], the latter class coincides with $ch_0(E)$. First we need a
more
explicit description of $ch^{inf}$ in terms of $Cyl^*(R,I)$.
We have a natural map $Cyl^*(R/I^n)\fib Cyl^*(R,I)_n$ and, by [3, 3.7]
also a natural map $Cyl^*(R,I)_{2(n^2+n-1)}\fib Cyl^*(R/I^n)$. 
Thus, if $\{u_n\}$ is a sequence
of elements $u_n\in G^1_n=Gl_s(\widetilde{Cyl}^1(R,I)_n)$ satisfying
$u_1=1$,
$u_n\faz^0_n(e_n)=\faz^1_n(e_n)u_n$, and compatible with the projections
$G^1_n@>>>G^1_{n-1}$, then the sequence $\{[\delta_n]\}$ of the $K_1$-classes
of the compatible elements:
$$
\delta_n=\faz^2_n(u_n)\faz^0_n(u_n)\faz^1_n(u_n)^{-1}e_n+1-e_n\tag{24}
$$
represents 
$ch^{inf}(E)$. Write $\Omega=\Omega(R)$,
$\faz^0=1+\sum_{i=1}^\infty D_i:R@>>>Cyl(R)\cong
\Omega/{\Omega^+}^\infty$ with $D_1=d$, $D_2=\phi$,  $D_i(R)\subset\Omega^i$.
Thus $\faz^0_n=1+\sum_{i=1}^{n-1}D_{i,n}:R/I^n@>>>Cyl^1(R,I)_n$ where
the $D_{i,n}$ are the induced maps. 
As in the proof of 4.2 above,  one checks that, for
$u_n=1+\sum_{i=1}^{n-1}\omega_{i,n}$, ($\omega_{i,n}\in\frac{\Omega^i}
{\Cal F^{n-i}\Omega^i})$, the condition that
$u_n\partial_n^0(e_n)=\partial_n^1(e_n)u_n$ is
equivalent to the sequence of identities:
$$
[e_n,\omega_{p,n}]=D_{p,n}(e_n)+\sum_{i=1}^{p-1}\omega_{i,n}D_{p-i,n}(e)
\qquad (1\le p\le n-1)\tag{25}
$$
Next one verifies that $\omega_{1,n}=(2e_n-1)de_n$ satisfies the identity
for $p=1$. Then one uses induction to show that the map 
$\psi_{p,n}:k@>>>\Omega^p/F_{n-p}\Omega^p$, 
$e\mapsto D_{p,n}(e_n)+\sum_{i=1}^{p-1}\omega_{i,n}D_{p-i,n}$ is a
derivation. It follows that if 
$\bar{\psi}_{p,n}:\Omega^1k@>>>\Omega^p/F_{n-p}\Omega^p$ is the induced
homomorphism, then $\omega_{p,n}=\bar{\psi}_{p,n}(\omega_{1,n})$ satisfies 
\thetag{25}. Hence $\{u_n=1+\sum_{i=1}^{n-1}\omega_{i,n}\}$ is a
compatible
sequence satisfying the requirements above. A formula for $\omega_{1,n}$
has
been given already; we have 
$\omega_{2,n}=\nabla(\omega_{1,n})$, where $\nabla(xdy)=x\phi(y)$. We shall
have no use for the terms of higher degree. From now on we fix $n$
and omit
it as a subscript. One uses \thetag{9} to calculate that the element
\thetag{24} is:
$$
\delta=1+e(de)^2\otimes(v_1v_2+v_1^2)+\dots\tag{26}
$$
Here and below, the \dots stand for higher degree terms. The isomorphism:
$$
H^2_{inf}(A,K_1)\cong H^2_{inf}(A,HC_0)
$$ 
of 4.0 is induced by the
logarithm
map $log:\ker(Gl(Cyl^2(R,I)_n)@>>>Gl(R/I^n))@>>>
\ker(HC_0(Cyl^2(R,I)_n)@>>>HC_0(R/I^n))$. Applying this map to \thetag{26}
we still get $log(\delta)=e(de)^2\otimes(v_1v_2+v_1^2)+\dots$ where
again the $\dots$ are higher degree terms, but not necessarily the same
as in \thetag{26}. The composite of the projection to the normalized
complex and the map $1\otimes p$ of the proof of 3.0 above applied
to $log(\delta)$ yields the element $e(de)^2$. Mapping the latter
through the periodicity map of 4.1 we obtain the element $2e$.\qed
\enddemo
\bigskip
\remark{Remark 4.5} Here is an explanation for the normalization factor
appearing in theorem 4.4. For simplicity we shall restrict our attention 
to the quasi-free case. Let $R$ be quasi-free, $E$ a finitely generated
right projective $\tilde{R}$-module and $\nabla$ a connection on $E$. 
Write $\nu=s'\nu'$ where $s':H_{dR}^2(R)\cong H^0_{dR}(R)\cong HP_0(R)$
is the periodicity isomorphism of 3.1. The proof above shows that:
$$
\nu'(ch^{inf}(E))=[tr(\nabla^2)]\in H^2_{dR}(R)
$$
Here $tr$ is the trace and $[\ \ ]$ denotes cohomology class.
Thus
the factor of $1/2$ needed to normalize $ch^{inf}(E)$ is the same as that
appearing
in Karoubi's definition of the Chern characters to de Rham cohomology
([11, 1.17]):
$$
ch_{2q}(E)=(1/2q!)[tr(\nabla^{2q})]\in H^{2q}_{dR}(R)
$$
\endremark
\bigskip
\bigskip
\subhead{5. The Jones-Goodwillie character}\endsubhead
\bigskip
The title of this section refers to the homomorphism:
$$
c_n:K_n(A)@>>>HN_n(A)\tag{27}
$$
going from $K$-theory to negative cyclic homology ($n\ge 1$).
The map \thetag{27} was defined for  arbitrary unital
rings in [13], where it was proved it induces an isomorphism at
the level of the rational relative groups associated to a nilpotent 
ideal. Here we shall restrict our attention to the particular case
when $A$ is a $\rat$-algebra. Because for algebras over $\rat$, 
relative $K$ and cyclic homology groups are 
$\rat$-spaces, (cf. [15]), it follows from [13] that if $A$ is a 
unital $\rat$-algebra and $I\subset A$ is nilpotent then:
$$
c_n:K_n(A,I)\iso HN_n(A,I)\tag{28}
$$
is an isomorphism.
It is easy to extend the map \thetag{27} to nonunital algebras and 
to show that \thetag{28} holds nonunitally also; see [3, 4.2] 
for details. 
Also in [3, 4.2] it is shown that \thetag{27} is the homomorphism
induced at the level of homotopy groups by a map of fibrant simplicial
sets:
$$
c:\zal_\infty\Cal NGl(A)@>>> SCN_{\ge 1}(A)\tag{29}
$$
Here $\zal_\infty$ is the Bousfield-Kan completion [1], $\Cal N Gl$ is the
nerve of the general linear group,
and $SCN_{\ge 1}$ is the simplicial abelian group the Dold-Kan correspondence
associates with the truncation of the negative cyclic chain complex. In
general,
if $(C_n,\delta)_{n\in\Bbb Z}$ is any chain complex, and $m\in\Bbb Z$,
we write $C_{\ge m}$ for the complex which is $C_n$ in each degree $n\ge
m$, $\delta C_m$ in degree $m-1$ and $0$ in degrees $<m-1$. The purpose of
this 
section is to prove the following:
\bigskip
\proclaim{Theorem 5.0} Let $A$ be a $\rat$-algebra. Then the 
Jones-Goodwillie map \thetag{29} fits into a homotopy fibration sequence 
of fibrant simplicial sets:
$$
\Bbb H_{inf}(A,\zal_\infty\Cal NGl)@>>>\zal_\infty\Cal NGl(A)@>c>>SCN_{\ge 1}(A)
$$
Here $\Hy_{inf}(A,\bgp)$ is the infinitesimal hypercohomology simplicial
set, as defined in 1.3 above.
\endproclaim 
\medskip
\proclaim{Corollary 5.1} Let $A=R/I$ be a presentation of $A$ as a 
quotient of a quasi-free algebra $R$. For $n\ge 2$, write 
$LK_n(A)=\pi_n(\underset{m}\to{\holi}\zal_\infty\Cal NGl(R/I^m))$.
Then there is an exact sequence:
$$\multline
\dots @>>> HN_{n+1}(A)@>>>LK_n(A) @>>>K_n(A)@>>>HN_n(A)
@>>>\dots 
 @>>>LK_2(A)@>>>\\
K_2(A)@>>>HN_2(A)@>>>
H_{inf}^0(A,K_1)@>>>K_1(A)@>>>HN_1(A)@>>>H_{inf}^1(A,K_1)
@>>>0\endmultline$$\endproclaim 
\demo{Proof of Corollary 5.1} The sequence of the corollary is the 
long exact sequence
of homotopy groups of the fibration of the theorem. Thus it suffices
to compute the hypercohomology groups $\Hy^n_{inf}(A,\bgp)$.
It follows from [3, 5.5] and from [1, IX.3.1] that 
$LK_n(A)=H_{inf}^n(A,\rat_\infty\Cal NE)$ where $E$ is the elementary 
group. By [1,4.4] we have a functorial homotopy fibration sequence:
$$
\zal_\infty\Cal NE @>>>\zal_\infty\Cal NGl @>>>\Cal N(K_1)
$$
By [1, 7.2], $\Hy_{inf}^n(A,\Cal N(K_1))=
H_{inf}^{1-n}(A,K_1)$
for $n=0,1$ and is zero otherwise.
On the other hand, because $\Hy_{inf}(A,-)$ is a $\holi$, it preserves 
homotopy fibration 
sequences of fibrant simplicial sets. The groups $\Hy^n_{inf}(A,\bgp)$
for $n\ge 0$ are now easily calculated from the long exact sequence of 
homotopy associated to the fibration sequence obtained from applying the 
functor $\Hy_{inf}(A,-)$ to the sequence above.\qed
\enddemo
\bigskip
\demo{Proof of Theorem 5.0} I claim that:
$$
\qquad\Hy_{inf}^n(A,SCN_{\ge1})=0 \qquad(n\le 0)\tag{30}
$$
Assume the claim holds. Let $F$ be the homotopy fiber of the map 
\thetag{29}. We know from [3, 5.2] that we have a map
$H_{inf}(-,X)@>>>X$ which is natural for functors going from rings to fibrant
simplicial sets, and that the functor $X\mapsto H_{inf}(-,X)$ preserves
homotopy fibration sequences. Using these facts, one obtains a commutative 
diagram where both rows are homotopy fibrations:
$$
\CD
\Hy_{inf}(A,F)@>>>\Hy_{inf}(A,\zal_\infty\Cal NGl)@>>>\Hy_{inf}(A,SCN_{\ge
1})\\
@VVV  @VVV  @VVV\\
F@>>>\zal_\infty\Cal NGl(A)@>>>SCN_{\ge 1}(A)\\
\endCD\tag{31}
$$
From the claim
we deduce that the first map in the top row of 
\thetag{31} is a homotopy equivalence. Next recall from [3] that for 
functors $X$ which map surjections
with nilpotent kernel into homotopy equivalences the map $H_{inf}(-,X)@>>>X$
is a weak equivalence. By \thetag{28}, this applies to $F$. It follows 
that the first vertical map of \thetag{31} is an equivalence. We have
shown that the claim implies the theorem. Next note that, by [3, 5.4.1]
 the claim
is equivalent to the assertion that the nonpositive hypercohomology groups
of the chain complex $CN_{\ge 1}$ are zero.
We shall actually show that $\Hy_{inf}^n(A,CN_\ge 1)=0$ for $n\le 0$,
and equals $HP_n(A)$ for $n\ge 1$. It is proven below (Corollary 5.3) that
$CN_{\ge 1}$ has the same hypercohomology as the complex having 
$\Omega^1_\natural$ in degree $1$, $b\Omega^1$ in degree zero, and zero
elsewhere. A standard spectral sequence argument shows that the hypercohomology
of this complex is $H_{inf}^*(A,HH_1)$. The latter is calculated by
means of the long exact sequence of cohomology groups associated with 
the following exact sequences of sheaves:
$$
0@>>>HH_1@>>>\Omega^1_\natural @>>>b\Omega^1@>>>0
\quad\text{and}\quad 0@>>>b\Omega^1@>>>O@>>>O/[O,O]@>>>0
$$
By Lemma 5.5 below and the first of the sequences above, we get that
\noindent $H^0_{inf}(A,HH_1)=0$ and
$H^n_{inf}(A,HH_1)=H^{n-1}(A,b\Omega^1)$.
From Theorem 2.3 and \thetag{18}, we obtain $H^*_{inf}(A,O)=0$;
from Theorem 3.0 and what we have just seen, we obtain that 
$H^1_{inf}(A,HH_1)=0$ and $H^n_{inf}(A,HH_1)=HP_n(A)$. This completes
the proof.\qed
\enddemo
\bigskip
\proclaim{Lemma 5.2} Let $R$ be a quasi-free $\rat$-algebra, 
$\Omega^*=\Omega^{*}(R)$, $\nabla$ as in \thetag{17}, 
$X=(\Omega^0\underset{b}\to{\overset{\natural d}
\to{\rightleftarrows}}\Omega^1_\natural)$
the $2$-periodic de Rham complex of [5], $CP^*=CP^*(R)$ the
$(b,B)$ periodic cyclic complex. Define
a homogeneous map $f:X@>>>CP$:
$$
f(\natural\omega)=\{(-\nabla B)^n(1-b\nabla)\omega\}_{n\ge 0},\quad
 f(a)=\{(-\nabla B)^n(a)\}_{n\ge 0}\quad
(\natural\omega\in\Omega^1_\natural, a\in\Omega^0.) $$ Then: \smallskip
\item{i)} The map $f$ is a well-defined supercomplex homomorphism.
\smallskip \item{ii)} For the natural projection $\pi:CP\fib X$, we have
$\pi f=1$. \smallskip \item{iii)} Put $h_m:\prod_{r\ge
0}\Omega^r@>>>\prod_{r\ge 0}\Omega^{r+2m+1}$, $$
h_m(\sum_{r=0}^\infty\omega_r)=
\sum_{r=1}^\infty(-1)^m\nabla(B\nabla)^m\omega_r\qquad{m\ge 0} $$ Then for
$h=\sum_{m=0}^\infty h_m$, we have $1-f\pi=
(B+b)h(1-f\pi)+h(1-f\pi)(B+b)$. \smallskip \item{iv)} For the filtration
$F$ of 1.2 above, we have: $$ \split
f_1(b\Omega^2+F_m\Omega^1)\subset\prod_{p\ge 0}
F_{m-(2p+1)}\Omega^{2p+1},\quad f_0(F_m\Omega^0)\subset\prod_{p\ge 0}
F_{m-2p}\Omega^{2p}\\ \text{\ \ and\ \ } h(1-f\pi)(\prod_{p\ge 0}
F_{m-p}\Omega^{p})\subset \prod_{p\ge 0} F_{m-1-p}\Omega^{p}.\\ \endsplit
$$ \endproclaim \demo{Proof} Part i) is straightforward; part ii) is
immediate. By 2.4, the map $h$ of iii) is a contracting homotopy for
$\ker\pi$; the assertion of iii) is immediate from this. Item iv) is
immediate from \thetag{12}. \enddemo \bigskip \proclaim{Corollary 5.3}
Let $A$ be any algebra. Then the infinitesimal hypercohomology of each of
the Hochschild, cyclic, negative cyclic and periodic cyclic chain
complexes associated to the functorial mixed complex
$\Omega^*\underset{b}\to{\overset{B} \to{\rightleftarrows}}\Omega^{*+1}$
equals the hypercohomology of the corresponding chain complex associated
to the functorial mixed complex $X$. 
\endproclaim 
\bigskip 
\remark{Remark 5.4} The same argument as in the proof of
Theorem 5.0 shows that, for any ring $A$
--not necessarily a $\rat$-algebra-- we have a fibration sequence:
$$
\Bbb H_{inf}(A,\rat_\infty\Cal NGl)@>>>\rat_\infty\Cal NGL@>c>>SCN_{\ge
1}(A\otimes \rat)
$$
Moreover the analogue of the exact sequence
of Corollary 5.1 holds, with the same proof.
The lower terms of this rational sequence give the exact sequence of [3,
6.2]), except that now we know the map $HN_1(A\otimes\rat)\fib
H^1(A,K_1\otimes\rat)$ is surjective. The surjectivity of this map was
conjectured in [3]; this true conjecture was derived from the wrong
conjecture that $H_{inf}^*(A, HN_1)=0$. In fact it is not hard to show,
using 5.3, that $\Bbb H_{inf}^*(A,HN_1)=\Bbb H_{inf}^*(A,HH_1)$ which
is as calculated in the proof of 5.0, and is nonzero in general.
\endremark
\bigskip
\proclaim{Lemma 5.5} Let $A$ be a $\rat$-algebra. Then: 
$$
H_{inf}^*(A,\Omega^1_\natural)=0 
$$ 
\endproclaim 
\demo{Proof} Let $A=R/I$
be a presentation of $A$ as a quotient of a quasi-free algebra. We have to
show that the cohomology of the complex $C(R,I,\Omega^1_\natural)$ is
trivial. By definition the latter complex is the inverse limit of the
complexes $\Omega^1_\natural Cyl(R,I)_n$ The proof has two main parts.  
The first part consists of showing that there is an isomorphism of
pro-vectorspaces: 
$$ \Omega^1_\natural Cyl^m(R,I)\cong\{\oplus_{r=1}^n
\frac{(\Omega^{r+1}\oplus\Omega^{r})\otimes T^r(V^m)}
{M^m_r+_r\negthickspace\Cal G^{m,n}}\}\tag{32} 
$$ 
Here $\Omega=\Omega(R)$,
$M^m_r=\{(b\omega\otimes x+u\otimes y+(-1) ^r \kappa(u)\otimes\zeta(y),
(-1)^{r+1}b(u)\otimes y): \omega\in\Omega^{r+1}, u\in \Omega^r, x,y\in
T^r(V^m)\}$, and:
$$
_r\Cal{G}^{m,n}=(\sum_{i_0+\dots
i_{r+1}+j\ge n-r}I^{i_0} dAI^{i_1}\dots I^{i_{r}}d(I^j)I^{i_{r+1}}\oplus
F^{n-r}\Omega^r)\otimes T^r(V^m)$$
Now to the first part. We consider
first the case when the ideal $I=0$. We have: 
$$ 
\split
\Omega^1(Cyl^m(R))&=\{\Omega^1(Cyl^m(R)_n)\}\\
&=\{\oplus_{r=0}^{n-1}\oplus_{p+q=r}\widetilde{\Omega}^p\otimes\Omega^q\otimes
T^r(V^m)\}\\ &=\{\Omega^1(\Omega_m)_{deg<n}\} 
\endsplit
$$
 Here
$\tilde{\Omega}^0=\tilde{R}$, $\tilde{\Omega}^p=\Omega^p$ if $p\ge 1$,
$\Omega_m=\Omega_m(R)=\oplus_{p\ge 0}\Omega^p\otimes T^p(V^m)$ is the
graded algebra of 1.1 above, and $_{deg<n}$ indicates a truncation of the
graded $\Omega_m$-module $\Omega^1(\Omega_m)$. Because the commutator
subspace is graded, we have 
$\{\Omega^1_\natural(Cyl^m(R)_n)\}\cong\{\Omega^1_\natural(\Omega_m)_{deg<n}\}$.
We write $\delta$ for the de Rham differential of $\Omega_m$.  Fix a
degree $r$. By repeated application of the Leibniz formula satisfied by
$\delta$, we see that every form of degree $r$ can be written as a sum of
of elements of the following types: $$ \align \text{Type} (p,0):&\quad
a_0da_1\dots da_p\delta a_{p+1}da_{p+2}\dots da_{r+1} \otimes x\\
\text{Type} (p,1):&\quad a_0da_1\dots da_{p-1}\delta(d a_{p})da_{p+1}\dots
da_{r}\otimes x\tag{33} \endalign $$ Here $x\in T^r(V^m)$, $a_i\in R$,
and if $\epsilon=0,1$, then for elements of type $(p,\epsilon)$, the index
$p$ runs between $\epsilon$ and $r$. One checks furthermore that the
sum of the subspaces generated by the elements of each type is direct,
i.e. that we have an isomorphism:
$$ 
(\oplus_{p=0}^{r}\Omega^{r+1}e_p\oplus
\oplus_{p=1}^r\Omega^re_p)\otimes
T^r(V^m)\cong_r\negthickspace\Omega^1(\Omega_m)\tag{34}
$$ 
Here we write
$\omega e_p$ for the element $\omega$ in the $p$-th summand; the
isomorphism maps $a_0da_1\dots da_{r+\epsilon}e_p$ to the element of type
$(p,\epsilon)$ of \thetag{33} ($\epsilon=0,1$), and
$_r\Omega^1(\Omega_m)$ is the homogenous summand of degree $r$. To compute
the degree $r$ part of $\Omega^1(\Omega)_\natural$, we must divide
$_r\Omega^1(\Omega_m)$ by the subspace generated by the elements of the
form $[a,\omega\otimes x]$ and $[da\otimes v,u\otimes y]$ for $a\in R$,
$\omega\in_r\negthinspace\Omega^1(\Omega_1)$,
$u\in_{r-1}\negthinspace\Omega^1(\Omega_1)$, $x\in
T^r(V^m)$, $v\in V^m$ and $y\in T^{r-1}(V^m)$. One calculates each of
these
commutators for each of the types of $\omega$ and $u$, and then pulls them
back through the isomorphism \thetag{34} to obtain that
$_r\Omega^1(\Omega_m)_\natural$ is isomorphic to the quotient of the left
hand side of \thetag{34} by the subspace $M'_r$ generated by the
following elements: 
$$ 
\align 
[\omega,a]e_p\otimes x,\quad &
([u,a]e_q+(-1)^{q+1}(udae_q+udae_{q-1}))\otimes x,\\ & \eta\otimes
xe_s-(-1)^{r-\epsilon}\kappa (\eta)\otimes\zeta(x)e_{s+1}\\
\endalign 
$$
Here $0\le p\le r$, $1\le q\le r$, $\epsilon\le s\le r-1$,
$\omega\in\Omega^{r+1}$, $u\in\Omega^r$, $\eta\in\Omega^{r+1-\epsilon}$,
$\epsilon=0,1$, $x\in T^r(V^m)$, and $\zeta$ is as in the proof of theorem
3.0 above. Thus $\eta\otimes xe_s\equiv (-1)^{r-\epsilon}\kappa
(\eta)\otimes\zeta(x)e_{s+1}$ $(\mod M'_r)$ whence every element in
$_r\Omega^1(\Omega_m)_\natural$ is the class of some element in
$(\Omega^re_r\oplus\Omega^{r+1}e_{r})\otimes T^r(V^m)$, or in other words
there is a surjective homomorphism going from the latter to the former.
One has to check that, upon identifying $\Omega^{r+1}\oplus\Omega^r$ with
$\Omega^{r+1}e_{r+1}\oplus\Omega^re_r$, we have
$M^m_r={M'}^m_r\cap((\Omega^{r+1}e_r\oplus\Omega^re_r)\otimes T^r(V^m))$.
The inclusion $\subset$ is immediate. To prove the other inclusion, let
$\omega^0e_r+\omega^1e_r\in((\Omega^{r+1}e_r\oplus \Omega^re_r)\otimes
T^r(V^m))\cap {M'}^m_r$. Then there exist elements
$u^0_i\in\Omega^{r+2}\otimes T^r(V^m)$, $u^1_j\in\Omega^{r+1}\otimes
T^r(V^m)$, $v^0_k\in\Omega^{r+1}\otimes T^r(V^m)$, and $v^1_l\in
\Omega^r\otimes T^r(V^m)$, $0\le i\le r$, $1\le j\le r$, $0\le k\le r-1$,
$1\le l\le r-1$, such that: 
$$ 
\omega^1e_r=\sum_{i=1}^r(b\otimes
1)(u_i^1)e_i+\sum_{i=1}^{r-1}(v_i^1e_i+
(-1)^r(\kappa\otimes\zeta)(v_i^1)e_{i+1})
$$ 
and:
$$
\omega^0e_r=\sum_{i=0}^r(b\otimes 1)(u^0_i)e_i+\sum_{i=1}^r(-1)^{i+1}
(u^1_ie_i+u^1_ie_{i-1})+\sum_{i=0}^{r-1}(v_i^0+(-1)^{r+1}(\kappa\otimes\zeta)
(v_i^0)) 
$$ 
In the first identity, for each $i\ne r$ the coefficient of
$e_i$ in the right hand side must be zero; thus if we apply
$((-1)^{r+1}\kappa\otimes\zeta)^{r-i}$ to it and then take the sum over
all $i$, we obtain the identity: 
$$ 
\omega^1=b\otimes
1(\sum_{i=1}^r((-1)^{r+1}\kappa\otimes\zeta)^{r-i}(u^1_i)) 
$$ 
By a similar
procedure, the second identity yields: 
$$
\omega^0=(b\otimes1)(\sum^r_{i=0}((-1)^r\kappa\otimes\zeta)^{r-i}(u_i^0)+
\sum_{i=1}^r(1+(-1)^r(\kappa\otimes\zeta))(-1)^{i+1}
((-1)^r\kappa\otimes\zeta)^{r-i}(u_i^1) 
$$ 
This proves the inclusion
$M'_r\cap(\Omega^{r+1}e_r\oplus \Omega^re_r)\subset M_r$. Thus the proof
of \thetag{32} for the case $I=0$ is complete. Now let $I\subset R$ be any
ideal. We have: 
$$ 
\Omega^1(Cyl^m(R,I)_n)=\oplus_{0\le
p,q<n}\frac{\Omega^p\delta\Omega^q}
{F^{n-p}\Omega^p\delta\Omega^q+\Omega^p\delta F^{n-q}\Omega^q}T^{p+q}(V^m)
$$ 
One checks that the filtrations 
$$
\split
\Omega^1(\Omega_m)\supset{\Cal
G''}^{m,n}&= \sum_{p,q\ge 0}(F^{n-p}\Omega^p\delta \Omega^q
+\Omega^p\delta
F^{n-q}\Omega^q) \otimes T^{p+q}(V^m)\\
\text{and}\qquad & {\Cal G'}^{m,n}=\sum_{n\le
p+q+i+j}(F^i\Omega^p\delta F^j \Omega^q)\otimes T^{p+q}(V^m)
\endsplit
$$
are
equivalent, whence $\Omega^1(Cyl^m(R,I))
\cong\frac{\Omega^1(\Omega_m)}{{\Cal G'}^{m,\infty}}$ as pro-vectorspaces.
For the degreewise identification \thetag{34}, we have: ${\Cal
G'}^{m,n}\cap (\Omega^re_s\otimes T^{r}(V^m))=F^{n-r}\Omega^re_s\otimes
T^{r}(V^m)$ and 
$$
\multline
\Cal{G'}^{m,n}\cap{\Omega^{r+1}e_s\otimes T^r(V^m)}=\\
\sum_{i_0+\dots i_{r+1}+j\ge n-r}I^{i_0}dAI^{i_1}\dots
I^{i_{s-1}}d(I^j)I^{i_s} \dots dA I^{i_{r+1}}
\endmultline
$$
We note that ${\Cal G'}^{m,n}$ is closed under $\kappa\otimes\zeta$, and
if
$u\otimes x\in{\Cal G'}^{m,n}\cap\Omega^{r+1}e_s$ $(s\ge 1)$, then
$b\otimes 1(u\otimes x)\in {\Cal G'}^{m,n}\cap(\Omega^{r}e_s\otimes
T^r(V^m))$. From
these observations, we see that the same argument as that given in the
proof of \thetag{32} for the case $I=0$ can be applied here, to prove now
that $(\Omega^{r+1}\oplus\Omega^re_r)\otimes T^r(V^m) \cap(M'_r+{\Cal
G'}^{m,n})=M_r+_r\negthickspace{\Cal G}^{m,n}$. This finishes the first
part of the proof. For the second part, consider the map $p$ of the proof
of theorem 3.0 above. Using the calculation of the first part of the
current proof and the arguments of the proof of 3.0, one deduces that
$1\otimes p$ induces a homotopy equivalence between the pro-complex
$\Omega^1_\natural(Cyl(R,I))$ and the pro-complex having
$\frac{\Omega^{m+1}\oplus\Omega^m}{N^m+\Cal G^{1,\infty}}$ in codimension
$m$, with $N^m=\{(b\omega+u-\kappa(u),(-1)^{m+1}b(u)):\omega
\in\Omega^{m+2}, u\in\Omega^{m+1}\}$ and with as coboundary map that
induced by $\Omega^{m+1}\oplus\Omega^m@>>>\Omega^{m+2}\oplus\Omega^{m+1}$,
$(\omega,u)\mapsto (d\omega,(-1)^m\omega+du)$. The latter complex is
contractible: a contracting homotopy is given by $(\omega,u)\mapsto
((-1)^{|\omega|}u,0)$. We see that this homotopy maps $N^*$ to $N^{*-1}$
and is continuous for the topology of the filtration $\Cal G^1$. This
concludes the second and last part of the proof.\qed 
\enddemo 
\bigskip

\bigskip
\subhead{Acknowledgements}\endsubhead 
Although most of the
research for this paper was carried out in Buenos Aires, the first
draft was written during a four week visit to M\"unster, where I was
a guest of Prof. J. Cuntz. I am grateful to him and to the people of
his team for their hospitality. My visit to M\"unster was
made possible by a special travel grant from the University of Buenos
Aires, which covered part of my transportation expenses and by 
SFB's 478 which supported my stay.
\linebreak
I am indebted to my colleague Prof.  
J.J. Guccione from whom I learned Lemma 2.4. 
Especial thanks to Dr. C. Valqui of M\"unster who made me aware of the
fact that the coefficients in the formula \thetag{21} receive the name of
Catalan numbers and appear in various parts of mathematics.

\enddocument